\documentclass[11pt]{amsart}
\usepackage{amsmath}
\usepackage{amssymb}
\usepackage{amsthm}
\usepackage{verbatim}
\makeatletter
\def\Large{\@setsize\Large{17\p@}\xivpt\@xivpt}
\def\LARGE{\@setsize\LARGE{20\p@}\xviipt\@xviipt}
\def\huge{\@setsize\huge{25\p@}\xxpt\@xxpt}
\def\Huge{\@setsize\Huge{30\p@}\xxvpt\@xxvpt}
\def\thebibliography#1{%
  \def\section{\@startsection{section}{1}{5.5mm}{0ex}{2ex}{\bf}}%
  \baselineskip 8pt
  \section*{\hfill References\hfill\@mkboth
 {REFERENCES}{REFERENCES}}\footnotesize\list
 {[\,\arabic{enumi}\,]}{\settowidth\labelwidth{[#1]}\leftmargin\labelwidth
 \advance\leftmargin\labelsep
 \usecounter{enumi}}
 \def\newblock{\hskip .11em plus .33em minus .07em}
 \sloppy\clubpenalty4000\widowpenalty4000
 \sfcode`\.=1000\relax}

\long\def\@caption#1[#2]#3{\par\addcontentsline{\csname
  ext@#1\endcsname}{#1}{\protect\numberline{\csname
  the#1\endcsname}{\ignorespaces #2}}\begingroup
    \@parboxrestore
    \normalsize
    \@makecaption{\csname fnum@#1\endcsname}{\ignorespaces #3}\par
    \vspace*{-2.5ex}\endgroup}
\makeatother
\newtheorem{thm}{Theorem}[section]
\newtheorem{prop}[thm]{Proposition}
\newtheorem{cor}[thm]{Corollary}
\newtheorem{lem}[thm]{Lemma}
\theoremstyle{remark}
\newtheorem{rem}[thm]{Remark}%
\renewcommand{\arraystretch}{1.6}
\newcommand{\bsharp}{\raisebox{.1ex}{$\boldsymbol\sharp$}}
\def\boplus{\,\mbox{\boldmath $ \oplus $}\,}
\def\botimes{\,\mbox{\boldmath $ \otimes $}\,}

\DeclareMathOperator{\GL}{GL}
\DeclareMathOperator{\PGL}{PGL}
\DeclareMathOperator{\End}{End}
\DeclareMathOperator{\Aut}{Aut}
\DeclareMathOperator{\Pic}{Pic}
\DeclareMathOperator{\Gal}{Gal}
\DeclareMathOperator{\tr}{tr}
\def\mathbb{\mathbf}
\newcommand{\Z}{{\mathbb Z}}            % the integers Z
\newcommand{\Q}{{\mathbb Q}}            % the rationals Q
\newcommand{\Qbar}{\overline{{\mathbb Q}}}      
\newcommand{\F}{{\mathbf F}}            % Used for finite fields
\newcommand{\m}{{\mathfrak{m}}}         % maximal ideal

\renewcommand{\O}{{\mathcal O}}  
\newcommand{\OX}{{\mathcal O}_X}  
\renewcommand{\L}{{\mathcal L}}  

\newcommand{\Xstar}{X_0^*(p)}           % X_0^*(p)
\newcommand{\modX}{{\mathcal X}_0(p)}   % a model for X_0(p)
\newcommand{\modXstar}{{\mathcal X}_0^*(p)}     % a model for X_0^*(p)
\newcommand{\redX}{\widetilde{\modX}}   % reduction of X_0(p)
\newcommand{\redXstar}{\widetilde{\modXstar}}   % reduction of X_0^*(p)
\newcommand{\T}{{\mathbf T}}            % Hecke algebra T
\newcommand{\wpbar}{{\widetilde{w}_p}}       % reduction mod p of w_p
\renewcommand{\P}{{\mathbf P}}          % Projective space P
\newcommand{\w}{{\omega}}               % omega
\newcommand{\f}{{\mathbf f}}            % finite part

\newcommand{\Y}{{\mathcal X}}
\newcommand{\A}{{\mathcal J}}

\newenvironment{pf}%
  {%\vspace{-2\theorempostskipamount}
\begin{list}%
  {Proof.}{
    \setlength{\topsep}{0mm}
    \setlength{\labelwidth}{-2mm}
    \setlength{\itemindent}{0mm}
    \setlength{\leftmargin}{0pt}
    \setlength{\listparindent}{5.5mm}
    }
  \item}{\qed\end{list}\par}
\def\qed{\hfill\penalty1\hspace*{0pt plus 1fill}$\blacksquare$\smallskip}
\begin{document}
\baselineskip 13pt
%
%\jot 0pt
%\baselineskip=2\baselineskip
%
%\noindent
%\underline{\makebox[\textwidth]{\bf Preliminary Version.\hfill\null}}
%\vspace{-4mm}
%
\vspace*{19mm}
\noindent
{\Large\bf Automorphisms of $ X_0^*(p) $}
\renewcommand{\thefootnote}{\fnsymbol{footnote}}
\footnote[0]{%
  14K02,14H25,%
  11G30,11G10,14H40,14K15. \\
  \indent
  This work was supported in part by
  Grant-in-Aid for JSPS Fellows 11-06151 and
  an NSF Postdoctoral Fellowship.}
\\[5mm]
{\sc Matthew Baker} and {\sc Yuji Hasegawa}
\par\vspace{10mm}
%\par\vspace{19mm}
\section{Introduction}
Let $ N $ be a positive integer, and let
$ X_0(N) $ be the modular curve corresponding
to the congruence subgroup
\[  \Gamma_0(N) =
   \bigg\{\bigg(\begin{matrix} a & b \\[-3.5mm] c & d \end{matrix}\bigg)
   \in {\mathrm{SL}}_2({\Z}) \biggm| c \equiv 0 \!\!\!\mod N \bigg\}. \]
It is known by \cite{at-le} that the full normalizer
$ {\mathfrak N}(\Gamma_0(N)) $ of $ \Gamma_0(N) $ in
$ \GL_2^+({\Q})=\{ \alpha\in\GL_2({\Q}) \mid \det\alpha > 0 \} $
contains the group $ \Gamma_0^*(N) $ generated by
all the matrices of the form
\[ \bigg(\begin{matrix} N'a & b \\[-3.0mm] Nc & N'd \end{matrix}\bigg),\quad
   a,b,c,d \in {\Z},\ N'|N,\ N'{}^2ad-Nbc=N', \]
and that the factor group $ W(N):=\Gamma_0^*(N)/\Gamma_0(N) $ is
an elementary 2-abelian group of order $ 2^{\omega(N)} $, %
where $ \omega(N) $ is the number of distinct prime divisors of $ N $. %
In case $ N $ is square-free, it is also known \cite{hel} that
$ \Gamma_0^*(N) $ is maximal as a Fuchsian group
(hence coincides with $ {\mathfrak N}(\Gamma_0(N)) $). %
Now put $ B(N)={\mathfrak N}(\Gamma_0(N))/\Gamma_0(N) $, %
which is naturally regarded as a subgroup of
the automorphism group $ \Aut X_0(N) $ of the curve $X_0(N)$.
\footnote{Throughout this paper, if $X$ is a curve defined over a
field $k$, we denote by $ \Aut X$ the group 
$\Aut_{k_s} ( X \times_k k_s )$ of automorphisms of $X_{k_s}$ over a
separable closure $k_s$ of $k$, where $X_{k_s} := X \times_k k_s$.  By 
Lemma~\ref{lem:DeligneMumford} below, this is the same as 
$\Aut_{K} ( X \times_k K) $ for any separably closed field extension $K$
of $k_s$.  Similarly, if $A$ is an abelian variety over $k$, then
$\End A$ will be shorthand for $\End_{k_s} ( A \times_k k_s )$. }

Recall that $ \Aut X_0(N) $ has been
completely determined \cite{ogg2}\cite{ke-mo}\cite{elk}; %
it agrees with $ B(N) $ whenever $ X_0(N) $ has genus $ g \geq 2 $ and
$ N \neq 37,63 $. %
(For these two exceptional cases, the group $ B(N) $ is
of index 2 in $ \Aut X_0(N) $.) %
In particular, if $ N $ is square-free, %
then we have $ \Aut X_0(N)=W(N) $
whenever $ X_0(N) $ has genus $ g \geq 2 $ and $ N \neq 37 $. %
\par
Let $ X_0^*(N) $ be the quotient curve of $ X_0(N) $ by $ W(N) $, %
and let $ g^*(N) $ be the genus of $ X_0^*(N) $. %
% Clearly $ X_0^*(N) $ corresponds to $ \Gamma_0^*(N) $. 
If $ N $ is a prime power, then $ X_0^*(N) $ coincides with
$ X_0^+(N) $, the quotient of $ X_0(N) $ by the Atkin--Lehner involution. %
We are interested in determining the group $ \Aut X_0^*(N) $
(for $ g^*(N) \geq 2 $). %
In the rest of this paper, we always assume $ N $ to be
{\it square-free}, unless otherwise specified. %
Then $ \Gamma_0^*(N) $ is a maximal Fuchsian group, %
so every nontrivial automorphism of $ X_0^*(N) $ is
necessarily exceptional in the sense that it
does not arise from a linear fractional transformation on
the complex upper half plane. %
Thus, if we assume that $ g^*(N) \geq 2 $, %
we expect the group $ \Aut X_0^*(N) $ to be very small. %
The purpose of this note is to prove that
when $ N $ is a prime number, this is indeed the case:
\begin{thm}\label{thm:main}
Let $ p $ be a prime such that $ g^*(p) \geq 2 $. %
Then
\begin{equation*}
\Aut X_0^*(p)=
\begin{cases}
{\Z}/2{\Z}\ & \text{if $ g^*(p)=2 $}; \\
\{ 1\}      & \text{if $ g^*(p)>2 $}.
\end{cases}
\end{equation*}
\end{thm}
We remark that Ogg's method in \cite{ogg2} for determining
$ \Aut X_0(p) $ does not readily generalize to $ X_0^*(p) $. %
Ogg's proof makes essential use of the action of
automorphisms on the cuspidal subgroup of $ J_0(p) $. %
However, as $X_0^*(p)$ has only one cusp, 
the cuspidal group in $ J_0^*(p) $ is trivial, %
so no information can be gained in this way. %
(Here $ J_0(p) $ and $ J_0^*(p) $ denote the Jacobians of
$ X_0(p) $ and $ X_0^*(p) $, respectively.)
\par
Our plan for proving Theorem 1.1 is as follows. %
In Section \ref{sec:str}, %
we investigate the structure of $ \Aut X_0^*(N) $ for $ N $ square-free. %
Using this, we show in Section \ref{sec:Aut=1 a.a.p} that
$ \Aut X_0^*(p) $ is trivial for almost all primes $ p $. %
Finally, by considering the reduction of $ X_0^*(p) $ modulo various primes, %
we show in Sections \ref{sec:red. mod l} and \ref{sec:red. mod p} that
$ \Aut X_0^*(p) $ is trivial for the remaining primes such that $ g^*(p)>2 $.

\medskip

To put Theorem~\ref{thm:main} in a larger context, we note that 
little is known in general about the problem of determining
the rational points on $ X_0^*(p) $. (For a discussion of what is 
currently known about this problem, including the existence of 
certain ``exceptional'' rational points for $p=191,311$, 
see S.~Galbraith's paper \cite{gal}).
On page 145 of his paper \cite{Mazur}, Mazur writes 
(adapting his terminology to ours): ``What further rational points
[in addition to the cusp and the rational CM points] does the curve
$X_0^*(p)$ possess?  This diophantine question (when the genus $g^*(p) > 0$)
is extremely interesting, since no known method appears to be applicable to
it, for {\it any} value of $p$.''
We hope that Theorem \ref{thm:main} will be useful
in future investigations of this problem.  In any case, if it had turned 
out that there were unexpected automorphisms of $X_0^*(p)$ for certain $p$, 
this might have ``explained'' the existence of certain 
exceptional rational points such as those found by Galbraith.

\medskip

Also, we remark that the curve $X_0^*(p^2)$ coincides with the curve which 
Mazur in \cite{Mazur2} calls $X_{\rm split}(p)$.  The rational points 
on this family of curves are also not known in general (see 
\cite{Momose} for some results in this direction).  
The problem of determining 
the rational points on $X_{\rm split}(p)$ is related to
a famous question formulated by Serre asking 
which elliptic curves $E/\Q$ have
the image of their mod $p$ Galois representations contained in a proper
subgroup of $\GL_2(\F_p)$ (in this case the normalizer of a split 
Cartan subgroup).  See \cite{Mazur2} for further details.  
It would be especially interesting for this reason to prove a result similar
to Theorem~\ref{thm:main} in which the automorphism group of $X_0^*(p^2)$ is
determined.  Our methods shed only partial light on this problem.

\medskip

{\small
The authors would like to thank Brian Conrad and the anonymous referee for
their useful comments on an earlier version of this manuscript.
}

\section{Structure of $ \protect\Aut X_0^*(N) $}\label{sec:str}
In this section, we investigate the structure of $ \Aut X_0^*(N) $
for $ N $ square-free.  It is known in this case (see \cite{DR}) that 
$J_0(N)$, and hence $J_0^*(N)$, is a semistable abelian variety over $\Q$.
We begin by recalling the following fact:
\begin{lem}\label{lem:AutX into EndJ}
Let $ X $ be a smooth, proper, and geometrically connected
algebraic curve of genus $ g\geq 2 $ over a field $ k $, 
and let $ J=\Pic^0_{X/k} $ be the Picard (Jacobian) variety of X over $k$.
Then $ \Aut_k (X) $ injects into $ \End_k(J) $.
\end{lem}
\begin{pf}
Suppose $ \phi\in\Aut_k(X) $ induces the identity map on $ J $. %
Then $ \phi $ acts trivially on the cotangent space at $ 0 $ in $ J $, %
which is canonically isomorphic to $ H^0(X,\Omega^1_X) $. %
It follows that $ \phi $ acts trivially on the image of $ X $
under the canonical map to projective space. %
Now if $ X $ is not hyperelliptic, then this canonical map is an embedding, %
which implies that $ \phi $ is identity on $ X $. %
If $ X $ is hyperelliptic, this argument shows that
$ \phi $ must be either the identity or the hyperelliptic involution $ h $. %
Since $ h $ acts as $ -1 $ on $ J $, %
it follows that the natural map from $\Aut_k(X)$ to $\End_k(J) $ is 
injective in the hyperelliptic case as well.
\end{pf}

We also have the following general fact:

\begin{lem}
\label{lem:DeligneMumford}
Let $ X $ as above be a smooth, proper, geometrically connected curve 
of genus $ g\geq 2 $ over a field $ k $, and let $K / k$ be an extension 
field.  Let $k_s$ denote the separable closure of $k$ in $K$.
Then every automorphism of $X_K$ is defined over $k_s$.
\end{lem}

\begin{pf}
Grothendieck's theory of the Hilbert scheme (see \cite{GrothDescent}) 
implies that there is a scheme 
${\rm \bf Aut}_{X/k}$, locally of finite type over $k$, 
which represents the functor associating 
to each $k$-scheme $S$ the set of $S$-isomorphisms from
$X \times_k S$ to itself.
It is enough to prove that ${\rm \bf Aut}_{X/k}$ is 
{\'e}tale over $k$, since this implies that every $K$-point comes from
a $k_s$-point.  This is a special case of 
Lemma~1.11 of \cite{DM}.
\end{pf}

\begin{cor}\label{cor:rationality}
Let $ X , k , J$ be as in Lemma~\ref{lem:AutX into EndJ}, 
and let $K / k$ be an extension field.
If $ \phi $ is a $K$-automorphism of $ X_K := X \times_k K $ such that
the induced $K$-endomorphism $ \phi^* $ of $ J \times_k K 
= \Pic^0_{X / K} $ is defined over $ k $, then $ \phi $ is defined over 
$ k $.
\end{cor}
\begin{pf}
We may clearly enlarge $K$ so that it is assumed to be separably closed.
If $k_s$ denotes the separable closure of $k$ in $K$, then by 
Lemma~\ref{lem:DeligneMumford}, we may assume
% \cite[Cor.~20.4]{AV} we have $\End_{k_s}(J) = \End_K(J)$,
% so we may assume without loss of generality that $K=k_s$.
without loss of generality that $K=k_s$.
Descent theory now tells us that $\phi$ is defined over $k$ if and only if
$\phi^g = \phi$ for all $g \in \Gal(k_s / k)$.  Since the natural map from
$\Aut_{k_s}(X \times_k k_s)$ to $\End_{k_s}(J\times_k k_s) $ is injective, 
it suffices to show that $(\phi^g)^* = \phi^*$ on $J \times_k k_s$ for all
$g$.  But $(\phi^g)^* = (\phi^*)^g$ by the definition of Picard functoriality,
so the hypothesis that $\phi^*$ is defined over $k$ gives us what we want.  

% The natural map $ \Aut X \to \End J $ is defined over $ k $. %
% So if $ \phi^\sigma \neq \phi $ for some 
% $ \sigma\in\Gal(\overline{k}/ k) $, % 
% then $ (\phi^*)^\sigma = (\phi^\sigma)^* $, %
% which by the lemma is different from $ \phi^* $. %
% This contradicts the fact that $ \phi^* $ is defined over $ k $.
\end{pf}
Lemma \ref{lem:AutX into EndJ} and Corollary~\ref{cor:rationality} are 
particularly useful when
the torsion part of $ \End(\Pic^0(X))^{\times} $ has a simple structure. %
For example, we obtain the following result:
\begin{prop}\label{prop:Aut X}
Let $ X $ be an algebraic curve over $ \Q $ of genus $ g \geq 1 $. %
Assume that its Jacobian $ J(X) $ is semistable over $ \Q $ and
that $ \End J(X) \botimes\Q $ 
is a product of totally real fields. %a
Then $ \Aut X $ is elementary $ 2 $-abelian$:$ %
$ \Aut X \cong
  {\Z}/2{\Z}\boplus\cdots\boplus{\Z}/2{\Z} $. %
Moreover, every automorphism of $ X $ is defined over $ \Q $.
\end{prop}
\begin{pf}
The first assertion follows from Lemma \ref{lem:AutX into EndJ}, %
since by assumption the torsion part of 
$ (\End J(X) )^{\times} $ is isomorphic to
$ {\Z}/2{\Z}\boplus\cdots\boplus{\Z}/2{\Z} $. %
The latter part is a consequence of Corollary \ref{cor:rationality}
and the following result of Ribet.
\end{pf}
\begin{thm}[Ribet \cite{rib1}]\label{thm:Ribet}
Let $ A $ be a semistable abelian variety over $ \Q $. %
Then every endomorphism of $ A $ is defined over $ \Q $.
\end{thm}

%
% Let $ J_0^*(N) $ be the Jacobian variety of $ X_0^*(N) $, %
% and let $ A $ be a $ \Q $-simple factor of $ J_0^*(N) $. %
% $ A $ has semistable reduction over $ \Q $, %
% so by Ribet's theorem every endomorphism of $ A $ is
% defined over $ \Q $. %
% It follows from \cite{shi},\cite{rib2} that $ \End A\botimes\Q $ is
% isomorphic to a totally real field. %
% Moreover, since $ N $ is assumed to be square-free, %
% the simple factors of $ J_0^*(N) $ are mutually non-isogenous.

Let $ J_0^*(N) $ be the Jacobian variety of $ X_0^*(N) $.  Since $N$ is 
square-free, $J_0^*(N)$ is semistable over $\Q$ (see \cite{DR}).  It
then follows from Theorem~\ref{thm:Ribet} and the results of 
\cite{rib2} that $ \End J_0^*(N) \botimes{\Q} $
is isomorphic to a product of totally real fields. %
Applying Proposition \ref{prop:Aut X}, we have
\begin{cor}\label{cor:Aut X_0^*(N)}
The group $ \Aut X_0^*(N) $ is elementary $ 2 $-abelian. %
Moreover, every automorphism of $ X_0^*(N) $ is defined over $ \Q $.
\end{cor}
As an application we can determine $ \Aut X_0^*(N) $
when $ J_0^*(N) $ is simple. %
Recall that $ X_0^*(N) $ is hyperelliptic
if and only if it is of genus $ 2 $ (cf.~\cite{ha-ha}). %
Therefore
\begin{cor}
If $ J_0^*(N) $ is simple, then
\begin{equation*}
\Aut X_0^*(N)=
\begin{cases}
{\Z}/2{\Z}      & \text{if $ g^*(N)=2 $}; \\
\{ 1 \}         & \text{if $ g^*(N)>2 $}.
\end{cases}
\end{equation*}
\end{cor}
Now let $ p $ be a prime number. %
Then $ X_0^*(p) $ is hyperelliptic
(i.e., $ g^*(p)=2 $) if and only if
\[ p=67,73,103,107,167,191. \]
Since $ J_0^*(p) $ is simple for these six primes
(see \cite[Table 5]{at-ting}), %
we conclude that Theorem \ref{thm:main} holds
when $ g^*(p)=2 $. %
\section{$ \protect\Aut X_0^*(p) $ is trivial for almost all $ p $}
\label{sec:Aut=1 a.a.p}
In the previous section we saw that $ \Aut X_0^*(N) $ is
an elementary $ 2 $-abelian group. %
The purpose of this section is to show that $ \Aut X_0^*(p) $ is trivial
for almost all primes $ p $. %

Throughout this section, let $p$ be a prime number and let $N$ be a 
square-free integer such that $g^*(p)$ and $g^*(N)$ are at least 2.

We need the following two lemmas.
\begin{lem}\label{lem:comm.}
Let $ u $ be an automorphism of $ X_0^*(N) $. %
Then, as endomorphisms of $ J_0^*(N) $, we have
\[ uT_l = T_l u \]
for each prime $ l $,
% such that $ l \nmid N $, 
where $ T_l $ is the $ l $-th Hecke operator.
\end{lem}
\begin{pf}
This follows from that fact that since $ N $ is assumed to be square-free,
$ \End J_0^*(N) $ is commutative (cf. \cite[Lem.\,2.6]{ke-mo}).
\end{pf}

Let $ \infty $ denote the cusp on $ X_0^*(N) $ which is the image of
the cusp $\infty$ on $X_0(N)$ under the natural map. %

\begin{lem}\label{lem:uc!=c}
If $ u $ is a nontrivial automorphism of $ X_0^*(N) $, %
then $ u\infty \neq \infty $.
\end{lem}
\begin{pf}
(See also \cite{ogg2}.)
The eigenvalues of $ u $ acting on
the space of holomorphic differentials of $ X_0^*(N) $ are all $ \pm 1 $. %
If they are all $ +1 $, %
then an argument using the canonical map to projective space
as in the proof of Lemma \ref{lem:AutX into EndJ} shows that
$ u=id $ on $ X_0^*(N) $. %
Similarly, if the eigenvalues are all $ -1 $, %
then $ X $ is hyperelliptic and $ u $ is the hyperelliptic involution. %
In this case, one can see explicitly using Weierstrass points 
that $ u\infty \neq \infty $. %
Recall that $X_0^*(N)$ is hyperelliptic exactly when $g^*(N)=2$.  
The Weierstrass points on a genus 2 curve $X$ are the fixed points
of the hyperelliptic involution, which is equivalent to saying that $P$
is a Weierstrass point if and only if there 
is a nonzero holomorphic differential vanishing to order 2 at $P$.
In our situation, by computing the $ q $-expansion at $ \infty $ of
each Hecke eigendifferential in the cases where $ g^*(N)=2 $, one finds that
$ \infty $ is not a Weierstrass point on $ X_0^*(N) $, so
it cannot be fixed by the hyperelliptic involution.
\par
Suppose, then, that $ u $ has eigenvalues $ +1 $ and $ -1 $. %
Then one can choose two differentials
$ \omega_1,\omega_2 \in H^0(X_0^*(N),\Omega^1) $ with
$ u^*\omega_1 = \omega_1, u^*\omega_2 = -\omega_2 $
which are normalized eigendifferentials
for the action of the Hecke algebra $\T$. 
The $q$-expansion at $\infty$ of each $\omega_i$ is of the form
\[
\omega_i = \Big(\sum_{m=1}^{\infty} a_m q^m\Big)\frac{dq}{q}
\]
with $a_1=1$.

We see that $\omega = \omega_1 + \omega_2$ does not vanish at $\infty$, but
its pullback $u^* \omega = \omega_1 - \omega_2$ {\it does} vanish at $\infty$.
Consequently, we must have $u\infty \neq \infty$.
\end{pf}

Note that since $\infty$ is the only cusp on $X_0^*(N)$ when $N$ is 
square-free, the statement
$u\infty \neq \infty$ is equivalent to the statement that 
$u\infty$ is not a cusp of $ X_0^*(N) $.

We now come to the following key lemma.
\begin{lem}\label{lem:key}
Let $ p $ be a prime number such that $ g^*(p) \geq 2 $, and
suppose there is a nontrivial automorphism $ u $ of $ X_0^*(p) $. %
Let $ \infty $ be the cusp of $ X_0^*(p) $, and
$ P $ a non-cuspidal geometric point of $ X_0^*(p) $. %
Then for each prime $ l $ such that $ l \neq p $, %
the divisor $ D_{l,P}:=(uT_l-T_l u)(\infty-P) $ on 
$ X_0^*(p)_{\Qbar} $ is non-zero but linearly equivalent to zero.
\end{lem}
\begin{pf}
The divisor $ D_{l,P} $ is linearly equivalent to zero
by Lemma \ref{lem:comm.}.  Also, note that by Lemma~\ref{lem:uc!=c}, 
$Q := u\infty$ is not a cusp.
Now suppose that the divisor $ D_{l,P} $ is identically zero, i.e., that 
$uT_l\infty + T_l uP = T_l u\infty + uT_l P$.  Since 
$T_l\infty = (l+1)\infty$, we see by applying $u$ to both sides that 
there can be no cancellation 
between points in the support of $uT_l\infty$ and $uT_l P$, because  
all points in the support of $T_l P$ are noncuspidal by the definition of
the correspondence $T_l$.
Therefore we must have complete cancellation between $uT_l\infty$ and
$T_l u\infty$.  In other words, we must have
$ (u T_l-T_l u)(\infty)=0 $, i.e.,  $ (l+1)Q=T_l(Q) $.
It follows that $ Q $ is a rational point of $ X_0^*(p) $
represented by a CM elliptic curve. %
By \cite[Thm.\,3.2]{gal}
(and the remarks following the proof of that theorem), %
we see that $ Q $ must be a Heegner point of class number one. 
(Note in particular that all noncuspidal rational points on $X_0^*(p)$ 
coming from rational points on $X_0(p)$ are Heegner points of 
class number one when $g_0^*(p) \geq 2$, and this happens only for
$p=67$ and $p=163$).
Let $ (E,C) $ represent a geometric point $P$ of $ X_0(p) $ 
which projects to $ Q $, so
that $Q$ itself can be thought of as corresponding to the set 
$\{ P,w_p P \}$ represented by the unordered pair 
$\{ (E,C), (E/C, E[p]/C) \}$. 
Since $ Q $ is a Heegner point of class number one, %
complex multiplication theory tells us that $ E/C \cong E $. %
If $ D $ is a cyclic subgroup of order $ l $ of $ E $,
then $ E/D $ is isomorphic either to $ E $ or to $ E/C $, %
since $ (l+1)Q=T_l(Q) $ as divisors on $ X_0^*(p) $. %
But $ E/C $ is isomorphic to $ E $, %
so we conclude that $ E/D \cong E $
for each of the $l+1$ cyclic subgroups $ D \subset E $ of order $ l $. %
The resulting degree $l$ maps from $E$ to $E / D \cong E$ give rise to
$l+1$ distinct elements of norm $l$ in $\End (E) / \{ {\rm units} \}$.
This is impossible, however.  To see this, let $\O = \End(E)$, let $K$ be
the fraction field of $\O$ (which is an imaginary quadratic field), 
and let $R$ be the ring of integers (maximal order) of $K$.
% Old, incorrect version: there are $ l+1$ such subgroups, 
% but strictly fewer than $l+1$ elements of norm $l$ in 
% $\End E / \{ \pm 1 \}$.
If the units of $\O$ coincide with those of $R$ (which is automatic except
when $K=\Q(i)$ or $\Q(\w)$, where $i$ and $\w$ are primitive fourth and 
third roots of unity, respectively), 
we would get $l+1$ distinct ideals of norm $l$ in $R$, which is clearly
impossible since $l+1 \geq 3$.  
We are left with the possibility that $R = \Z[i]$ or $\Z[\w]$ and
$\O$ is a nonmaximal order in $K$.  However, it is easy to see that 
these exceptional cases are also impossible using the fact that
in $\Z[i]$, 2 ramifies and 3 is inert, and in $\Z[\w]$, 3 ramifies and 
$2$ and $5$ are inert.
\end{pf}
\begin{cor}\label{cor:X0*(p):6-gon}
Suppose there is a nontrivial automorphism $ u $ of $ X_0^*(p) $. %
Then there is a finite morphism $ f{:}\,X_0^*(p) \to {\P}^1 $
of degree at most $ 6 $, defined over $ \Q $. %
Moreover, every nontrivial automorphism of $ X_0^*(p) $ has
at most $ 12 $ fixed points.
\end{cor}
\begin{pf}
We let $ l=2 $ in the above lemma. %
For the first assertion, take $ P=u\infty $. %
For the latter, take $ P \neq \infty, u\infty $ and consider
a morphism $ f{:}\,X_0^*(p) \to {\P}^1 $ such that $ (f)=D_{l,P} $. %
Then $ u^*f \neq f $, since a comparison of divisors shows that 
$u^*f$ has a zero at $\infty$ whereas $f$ does not.
Therefore the assertion follows from the following general fact.
\end{pf}
\begin{lem}\label{lem:2d>=r}
Let $ X $ be a smooth projective curve over an algebraically closed
field $k$, 
and let $ u $ be a nontrivial automorphism of $ X $
having $ r $ fixed points. %
If there is a finite morphism $ f{:}\,X \to {\P}^1 $
of degree $ d $ such that $ u^*f \neq f $, %
then $ r \leq 2d $.
\end{lem}
\begin{pf}
%% (See also \cite[p.\,261, Prop.\,V.1.4]{fa-kr}.)
Let us write $ g=u^*f-f $. %
Assume first that $ g $ is constant. %
Then all the fixed points of $ u $ must be poles of $ f $, %
so $ 2d \geq d \geq r $. %
Next assume that $ g $ is non-constant. %
Then $ 2d-r' \geq \deg g $, where $ r' $ is
the number of fixed points of $ u $ which are poles of $ f $. %
Furthermore, since every fixed point of $ u $
which is {\it not} a pole of $ f $
is a zero of $ g $, we must have $ 2d-r' \geq r-r' $.
\end{pf}
Now let $ l $ be a prime such that $ l \neq p $.  By abuse of notation, we
let $X_0(p)$ denote the standard model of $X_0(p)_\Q$ over $\Z [1/p]$, 
and we let $ \widetilde{X}_0(p) $ be the mod $l$ reduction of 
$ X_0(p) $, which is again a smooth curve.

By Corollary \ref{cor:X0*(p):6-gon}
there must exist a nonconstant morphism
$ X_0(p) \to {\P}^1 $ of degree at most 12 over $ \Q $
whenever $ \Aut X_0^*(p) $ is nontrivial. %

We claim that there necessarily exists a nonconstant morphism of 
degree at most 12 from $\widetilde{X}_0(p)$ to $\P^1_{\F_l}$ defined over
$\F_l$.  This follows from the following general lemma (c.f.~ 
\cite[Chapter III, Lemma 2.6]{Mazur}):

\begin{lem}
\label{semicont}
Let $R$ be a discrete valuation ring with field of fractions $K$, 
residue field $k$, and uniformizing parameter $\pi$.
Let $X$ be a regular scheme which is proper and 
flat over ${\rm Spec}(R)$ whose
generic fiber $X_K$ is a smooth and geometrically connected 
curve over $K$.  Let $X_k$ denote
the (possibly singular) closed fiber of $X$.
Then any line bundle $\L_K$ of degree $d$
on $X_K$ extends to a line bundle $\L$ on $X$, and the pullback $\L_k$ to
$X_k$ is a degree $d$ line bundle with $\dim_K H^0(X_K,\L_K) \leq \dim_k
H^0(X_k,\L_k)$.
\end{lem}

\begin{pf}
With our hypotheses, it follows from 
\cite[II, 6.11]{Hartshorne} that there is a natural isomorphism
between the  divisor class group of $X$ and the group Pic($X$).  
The surjectivity of the natural map ${\rm Pic}(X)\to {\rm Pic}(X_K)$
then follows from \cite[II, 6.5(a)]{Hartshorne}.  (If $X/R$ is smooth
then this map is also injective).
Let $\L$ be a line bundle on $X$ extending the degree $d$ line bundle $\L_K$
on $X_K$.
If we define the degree of a line bundle $\L$ to be
$\chi(\L)-\chi(\OX)$, then this degree is constant on fibers 
because $X/R$ is flat; by Riemann--Roch this degree must be $d$.

We have $\dim_K H^0(X_K,\L_K) = \dim_K H^0(X,\L)\otimes_R K$ because
cohomology commutes with the flat base extension $R\to K$ by 
\cite[III, 9.3]{Hartshorne}.

The pullback (restriction) $\L_k$ of $\L$ is a line bundle of degree $d$ on
$X_k$, and we have $\dim_K H^0(X,\L)\otimes_R K \leq \dim_k
H^0(X_k,\L_k)$ by the semicontinuity theorem 
(c.f.~\cite[III, 12.8]{Hartshorne}).

\end{pf}

\begin{cor}
Let $X$ be as in the above lemma, and assume furthermore that
the closed fiber $X_k$ is smooth.  Suppose $X_K$ admits a morphism of
degree at most $d$ to $\P^1$ defined over $K$.  Then $X_k$ admits
a morphism of degree at most $d$ to $\P^1$ defined over $k$.
\end{cor}

\begin{pf}
This follows directly from the lemma, since if $Y$ is any
smooth, proper, geometrically connected curve over a field $k$, 
then there exists a 
morphism of degree at most $d$ from $Y$ to $\P^1$ defined
over $k$ iff there exists a line bundle $\L$ of degree at most $d$ on
$Y$ such that $\dim_k H^0(Y,\L) \geq 2$.
\end{pf}

In our situation, since $\bsharp \P^1(\F_{l^2}) = l^2 + 1$, we see that 
the number $ \bsharp\widetilde{X}_0(p)({\F}_{l^2}) $
cannot exceed $ 12(l^2+1) $. %
On the other hand, Ogg \cite{ogg1} gives the following estimate:
\[ \bsharp\widetilde{X}_0(p)({\F}_{l^2}) \geq
   \frac{l-1}{12}(p+1)+2. \]
Therefore, if $ \Aut X_0^*(p) $ is nontrivial, %
we must have the inequality
\[ 12(l^2+1) \geq \dfrac{l-1}{12}(p+1)+2. \]
Setting $ l=2 $, we obtain the following:
\begin{cor}
If $ p > 695 $, then $ \Aut X_0^*(p) $ is trivial.
\end{cor}
The next result follows from Corollary \ref{cor:X0*(p):6-gon}:
\begin{cor}\label{cor:num>12}
If $ J_0^*(p) $ has a $\Q$-simple factor of
dimension larger than $ (g^*(p)+5)/2 $, %
then $ \Aut X_0^*(p) $ is trivial.
\end{cor}

\begin{pf}
Suppose $u$ is a nontrivial automorphism of $X_0^*(p)$, which we know must
have order 2.  By Corollary \ref{cor:X0*(p):6-gon}, the number $r$ of 
fixed points of $u$ is at most 12.  Let $X'$ be the quotient of 
$X_0^*(p)$ by $u$, let $g'$ be its genus, and let $J'$ be the Jacobian of
$X'$.  Applying the Riemann--Hurwitz formula
to the degree 2 map $X_0^*(p) \to X'$, we find that
$2g' - 1 \leq g^*(p) \leq 2g' + 5$, so that 
\begin{equation}
\label{genus-inequalities}
(g^*(p) - 5)/2 \leq g' \leq (g^*(p) + 1)/2.
\end{equation}
In particular, if $m > (g^*(p) + 1)/2$ and 
if $ J_0^*(p) $ has a $\Q$-simple factor $A$ of
dimension equal to $m$, then $g^*(p) = \dim J_0^*(p) \geq m + 
\dim J' \geq m + (g^*(p) - 5)/2$, which implies that 
$m \leq (g^*(p) + 5)/2$ as desired. 
\end{pf}
\begin{rem}
The $ \Q $-simple factors of $ J_0^*(N) $
are the same as the $\Qbar$-simple factors of $ J_0^*(N) $
whenever $ N $ is square-free, since in this case
$ J_0^*(N) $ is semistable over $ \Q $
(c.f.~\cite{rib1}; see also Theorem \ref{thm:Ribet}). %
\end{rem}
\begin{rem}
Let $ N $ be a positive integer. %
We may compute the $ \Q $-simple splitting of $ J_0^*(N) $
by factoring characteristic polynomials of Hecke operators
acting on the space $ S_2^*(N) $ of cuspforms of weight 2
on $ \Gamma_0^*(N) $. %
For instance, suppose we are given an element $ T $ of
the Hecke ring $ {\Z}[\{T_l\}_{l\nmid N}] $
whose action on $ S_2^*(N) $ has
{\it square-free} characteristic polynomial $ \Psi_T(x) $. %
Let $ \Psi_T=\Psi_1\cdot \Psi_2\cdots\Psi_s $ be
the factorization of $ \Psi_T $ over $ \Q $. %
Then
\[ \deg\Psi_1+\deg\Psi_2+\cdots+\deg\Psi_s \]
gives the $ \Q $-simple splitting of $ J_0^*(N) $. %
\par
Now consider the case $ g^*(p)>2 $. %
(For a discussion of how to compute $g^*(p)$, see \cite[p.~20]{Ogg4}).
Using the above observation, %
we find that there are 39 primes in the range $ p < 695 $
for which $ J_0^*(p) $ is not simple
(see Table \ref{tab:simple splitting}). %
Applying Corollary \ref{cor:num>12} to these cases, %
we see that $ \Aut X_0^*(p) $ is trivial except possibly for
\begin{align*}
p={} & 163,193,197,211,223,227,229,269,331,347, \\
     & 359,383,389,431,461,563,571,607.
\end{align*}
\end{rem}
\begin{table}[h]
\caption{Simple splitting of $ J_0^*(p) $}
\label{tab:simple splitting}
{\renewcommand{\arraystretch}{1.10}
\begin{center}
\begin{tabular}{|r|p{19mm}|}
\noalign{\hrule height 0.5pt}
$ p $ \hfil & Splitting \\
\noalign{\hrule height 0.5pt}
163 & $1+5$ \\
193 & $2+5$ \\
197 & $1+5$ \\
211 & $3+3$ \\
223 & $2+4$ \\
227 & $2+3$ \\
229 & $1+6$ \\
269 & $1+5$ \\
277 & $1+9$ \\
331 & $1+3+7$ \\
\noalign{\hrule height 0.5pt}
\end{tabular}
\hspace{0.1mm}
\begin{tabular}{|r|p{19mm}|}
\noalign{\hrule height 0.5pt}
$ p $ \hfil & Splitting \\
\noalign{\hrule height 0.5pt}
347 & $1+2+7$ \\
359 & $1+1+4$ \\
373 & $1+12$ \\
383 & $2+6$ \\
389 & $2+3+6$ \\
397 & $2+13$ \\
431 & $1+3+4$ \\
439 & $2+9$ \\
443 & $1+1+12$ \\
457 & $2+15$ \\
\noalign{\hrule height 0.5pt}
\end{tabular}
\hspace{0.1mm}
\begin{tabular}{|r|p{19mm}|}
\noalign{\hrule height 0.5pt}
$ p $ \hfil & Splitting \\
\noalign{\hrule height 0.5pt}
461 & $2+3+7$ \\
467 & $1+12$ \\
491 & $2+10$ \\
499 & $2+16$ \\
503 & $1+10$ \\
523 & $2+15$ \\
547 & $2+18$ \\
557 & $1+18$ \\
563 & $3+3+9$ \\
571 & $3+6+10$ \\
\noalign{\hrule height 0.5pt}
\end{tabular}
\hspace{0.1mm}
\begin{tabular}{|r|p{19mm}|}
\noalign{\hrule height 0.5pt}
$ p $ \hfil & Splitting \\
\noalign{\hrule height 0.5pt}
587 & $5+13$ \\
593 & $1+18$ \\
599 & $2+11$ \\
607 & $5+7+7$ \\
613 & $5+18$ \\
647 & $2+14$ \\
653 & $7+17$ \\
659 & $1+16$ \\
677 & $1+2+18$ \\
    &        \\
\noalign{\hrule height 0.5pt}
\end{tabular}
\end{center}
}
\end{table}
\section{Reduction of $ X_0^*(p) $ modulo $ l $}
\label{sec:red. mod l}
In this section we consider the reduction of
$ X_0^*(p) $ and its reduction modulo $ l $, where
$ l $ is a prime such that $ l \neq p $. %
Let $ S_2^*(p) $ be the space of
cuspforms of weight 2 on $ \Gamma_0^*(p) $. %
For each Hecke-stable subspace $ S $ of $ S_2^*(p) $, define
\[ \nu_{S,l^n}=1+l^n-\tr T_{l^n}|S+
          \begin{cases}
            \,l\cdot\tr T_{l^{n-2}}|S\ %
                & \text{if $ n \geq 2;$} \\
            \,0 & \text{otherwise}.
          \end{cases}
\]
Formulas for computing the traces of Hecke operators are
given in \cite{hij},\cite{yam}.
\medskip

We note that since the $\Q$-simple factors of $J_0(p)$ (and hence of
$J_0^*(p)$) are Hecke-stable and appear with multiplicity one (see 
\cite[Chapter II, Section 10]{Mazur}), it follows
that if $Y$ is any curve of positive genus covered by $X_0^*(p)$, then
the image of the Jacobian of $Y$ in $J_0^*(p)$ via Picard 
functoriality is Hecke-stable.  This implies that the subspace $S$ of
$S_2^*(p)$ (which is canonically identified with the cotangent space of 
$J_0^*(p)$ at the origin) corresponding to the cotangent space of 
$J$ at the origin is Hecke-stable.  We say colloquially in this 
situation that the Jacobian of $Y$ ``comes from'' $S$.

\medskip

Before giving the next proposition, we recall the following fact, which
is proved, for example, in \cite{Youssefi}.
Suppose $R$ is a discrete valuation ring with fraction field $K$, and
that $f : X \to Y$ is a finite morphism of projective, smooth, geometrically 
connected curves over $K$.  Assume furthermore that $X(K)$ is nonempty.  
Then if $X$ has good reduction, $Y$ does as well.  (We say a curve
$X$ as above has {\it good reduction} if there exists a scheme
of finite type which is proper and smooth over $R$ and 
whose generic fiber is $X$.)

\begin{prop}
Suppose there is a nontrivial automorphism $ u $ of $ X_0^*(p) $. 
(By the remarks at the end of section 2, we may assume that $g^*(p) > 2$,
so that $X_0^*(p)$ is not hyperelliptic.)
Let $ Y=X_0^*(p)/\langle u \rangle $ be
the quotient of $ X_0^*(p) $ by the action of $ u $, which is a curve of
genus $g \geq 1$, and let $ \tilde{Y} $ be the reduction of $ Y $ 
modulo $ l $. %
If the Jacobian $J$ of $ Y $ comes from $ S $, %
we must have
$ \bsharp\tilde{Y}({\F}_{l^n})=\nu_{S,l^n} $.
\end{prop}
\begin{pf}
% Let $ g $ be the genus of $ Y $. 
We can write
$ \bsharp\tilde{Y}({\F}_{l^n})=
  1+l^n-\sum_{i=1}^g(\alpha_i{}^n+\overline{\alpha}_i{}^n) $, where
the $ (\alpha_i,\overline{\alpha}_i) $ are conjugate pairs of
eigenvalues of Frobenius acting on the $l'$-adic Tate module of $J$
for any prime $l' \neq l$. %
On the other hand, since $ u $ is defined over $ \Q $, %
the Jacobian of $ Y $ may be considered, up to $ \Q $-isogeny, as
an abelian subvariety over $ \Q $ of $ J_0^*(p) $. %
Thus from the Eichler--Shimura congruence relation we conclude that
the set $ \{\alpha_i+\overline{\alpha}_i\}_{i=1}^g $ gives
the eigenvalues of $ T_l|S $. %
Using the standard recurrence relations for the Hecke eigenvalues, %
we obtain the desired equality.
\end{pf}
It follows from this proposition that
$ S $ cannot correspond to the Jacobian of $ Y $ once we know
$ \nu_{l^n}:=\bsharp\tilde{X}_0^*(p)({\F}_{l^n})>2\cdot\nu_{S,l^n} $
for some prime power $ l^n $. %
Checking this inequality for all $ S \neq \{1\} $ such that
$ (g^*(p)-5)/2 \leq g_S \leq (g^*(p)+1)/2 $ %
(see (\ref{genus-inequalities}) above),
where $ g_S $ is the dimension of $ S $, %
we see that $ \Aut X_0^*(p) $ is trivial for
$ p=163,193,197,269,277,331,347,359,383,461 $, $ 563,607 $
(Table \ref{tab:X->Y}).
% ; see also Table \ref{tab:simple splitting}
% for subspaces $ S \neq \{ 1\} $ satisfying
% $ (g^*(p)-5)/2 \leq g_S \leq (g^*(p)+1)/2 $).
%
%
%
\begin{table}[h]
\caption{}
\label{tab:X->Y}
{\renewcommand{\arraystretch}{1.10}
\begin{center}
\begin{tabular}{|c|c|c|c|c|}
\noalign{\hrule height 0.5pt}
 $p$ & $l^n$ & $\nu_{l^n}$ & $\nu_{S,l^n}$ & $g_S$ \\
\noalign{\hrule height 0.3pt}
163 & 2 &  8 & 3 & 1 \\
193 & 5 & 14 & 6 & 2 \\
197 & 3 & 12 & 4 & 1 \\
269 & 3 &  9 & 4 & 1 \\
331 & 4 & 20 & 5 & 3 \\
\noalign{\hrule height 0.5pt}
\end{tabular}
\hspace{0.1mm}
\begin{tabular}{|c|c|c|c|c|}
\noalign{\hrule height 0.5pt}
 $p$ & $l^n$ & $\nu_{l^n}$ & $\nu_{S,l^n}$ & $g_S$ \\
\noalign{\hrule height 0.3pt}
331 & 4 & 20 & 8 & 4 \\
347 & 3 & 11 & 4 & 3 \\
359 & 4 & 17 & 8 & 1 \\
359 & 4 & 17 & 8 & 1 \\
359 & 5 & 10 & 4 & 2 \\
\noalign{\hrule height 0.5pt}
\end{tabular}
\hspace{0.1mm}
\begin{tabular}{|c|c|c|c|c|}
\noalign{\hrule height 0.5pt}
 $p$ & $l^n$ & $\nu_{l^n}$ & $\nu_{S,l^n}$ & $g_S$ \\
\noalign{\hrule height 0.3pt}
383 & 9 & 38 & 15 & 2 \\
461 & 9 & 42 & 19 & 5 \\
563 & 4 & 31 & 15 & 6 \\
607 & 4 & 27 & 11 & 7 \\
607 & 4 & 27 & 12 & 7 \\
\noalign{\hrule height 0.5pt}
\end{tabular}
\end{center}
}
\end{table}
\par
It remains to determine Aut $X_0^*(p)$ for the seven primes
\[
p = 211,223,227,229,389,431,571.
\]
\section{Reduction of $ X_0^*(p) $ modulo $ p $}
\label{sec:red. mod p}
In this section, we consider the reduction of $\Xstar$ modulo $p$.

\medskip

Let $\modX$ be the Deligne--Rapoport model of $X_0(p)$ over $\Z_{(p)}$, 
whose special fiber
$\redX$ can be identified with a union of two copies of 
$X_0(1)_{\F_p} \cong \P^1_{\F_p}$ intersecting transversely at the 
supersingular geometric points.  
The automorphism $w_p$ of $\modX$ induces an automorphism $\wpbar$ 
of the special
fiber $\redX$ which interchanges the components $Y_1$ and $Y_2$ of
$\redX$.  
It takes each supersingular geometric point to its $\F_{p^2}$-conjugate, and 
in particular
the fixed points of $\wpbar$ correspond precisely to the supersingular 
$j$-values in characteristic $p$ which are defined over $\F_p$
(See \cite{Ogg3} and \cite{DR}).

\medskip

The quotient of $\modX$ by $w_p$ is an arithmetic surface $\modXstar$ 
whose generic fiber is the smooth curve $\Xstar$, and whose closed fiber
$\redXstar$ is a 
copy of the $j$-line that intersects itself transversely
at the $g^*=g^*(p)$ geometric points corresponding to conjugate pairs of 
supersingular elliptic curves defined over $\F_{p^2}$ but not $\F_p$.
% (see Appendix C of \cite{ha-ha} for another discussion of this).  
The fact that the closed fiber of $\modXstar = \modX / w_p$ 
is the quotient of 
$\redX$ by the action of $w_p$ follows from the discussion in 
\cite[Appendix A7]{KM}, especially Prop.~A7.1.3.\footnote{Specifically, 
we are using the fact that (a) the order 2 of the group 
generated by $w_p$ is invertible in $\F_p$ since $p\neq 2$, and
(b) the order of all stabilizer groups of all geometric points of 
$X_0(p)$ are invertible in $\F_p$ since $p\neq 2,3$.  
The reason we need (b) is that the Deligne--Rapoport model is only a 
coarse moduli scheme, but one can show under the hypothesis in (b) that 
formation of the coarse moduli scheme associated to an algebraic stack 
commutes with the base change in question.}

One can see that the fixed points of $\wpbar$ on $\redX$ get mapped 
to smooth points of $\redXstar$ by using the fact that the components
$Y_1$ and $Y_2$ of $\redX$ intersect transversely, plus the following local
calculation: if $k$ is any field, the subring of $k[[x,y]]/(xy)$ consisting
of all elements left invariant by the automorphism which interchanges $x$ and
$y$ is isomorphic to $k[[t]]$, a power series ring in one variable.

\medskip

We assume for the rest of this section that $g^*(p) > 0$. 

\medskip

Note that if $p \not\equiv 1$ mod $12$, then $\modX$ is not regular:
the supersingular points corresponding to elliptic curves with extra
automorphisms, i.e., to elliptic curves with $j$-invariant 0 or 1728, are
not regular.  However, the quotient of $\modX$ by $w_p$ behaves better in
this respect:

\begin{lem}
The arithmetic surface $\modXstar$ is the minimal proper regular model 
over $\Z_{(p)}$ for $\Xstar$.
\end{lem}

\begin{pf}
We first claim that $\modXstar$ is indeed regular.
Note that $\modXstar$ is normal and its closed fiber is generically 
smooth.  Since the generic fiber $\Xstar$ is smooth, the only points which
are potentially not regular are the finitely many singular points 
on the closed fiber.
Let $s\in\redXstar(\F_{p^2})$ be such a point, corresponding to a 
supersingular $j$-invariant $j(s)$ defined over $\F_{p^2}$ but not $\F_p$.
In particular, $j(s)$ is not 0 or 1728, and therefore $\pi^{-1}(s)$ 
consists of two distinct {\it regular} points $s_1,s_2$
of $\modX$, where $\pi$ is the canonical map from $\modX$ to $\modXstar$.
Since $s_1$ and $s_2$ are not fixed points of $w_p$,
it follows from \cite[Thm.\,A7.1.1]{KM} that $\pi$ is {\'e}tale over $s$.
Therefore $s$, being the image of a regular point by an {\'e}tale map, is
regular. 

To see that $\modXstar$ is minimal, it suffices by Castelnuovo's criterion 
(\cite[Thm.\,3.1]{Chinberg})
to note that $\modXstar$ has no exceptional divisors in the sense of 
\cite[Def.\,1.3]{Chinberg}.  This follows immediately
from the fact that the special fiber is irreducible, and therefore by
\cite[Cor.\,4.1]{Chinberg} has 
self-intersection zero.  
\end{pf}

Let $u$ be an automorphism of $\Xstar$.  By the defining property of
the fact that $\modXstar$ is a minimal proper regular model 
(see \cite[Def.~1.1]{Chinberg}), 
$u$ extends to an automorphism of the scheme $\modXstar$.
% preserving the 
% largest subscheme $(\modXstar)^0$ of $\modXstar$ smooth over $\Z_{(p)}$.  
It therefore induces an
automorphism $\widetilde{u}$ of the special fiber $\redXstar$ and 
an automorphism $\alpha$ of
the normalization $\P^1_{\F_p}$ of $\redXstar$ which permutes
the set of $2g^*$ supersingular geometric points defined over $\F_{p^2}$ 
but not $\F_p$.

\medskip

We claim that if $u$ is not the identity map on the generic fiber $\Xstar$,
then $\widetilde{u}$ is not the identity on the special fiber.  

This follows from Lemma~\ref{lem:AutX into EndJ} and the fact that 
$J_0^*(p)$ has semistable reduction at $p$, together with the following
two basic results:

\begin{thm}[Raynaud]
Let $R$ be a discrete valuation ring with fraction field $K$.  
Let $\Y/ R$ be a proper, flat, and regular curve whose 
generic fiber is smooth and geometrically irreducible and whose closed fiber 
is geometrically reduced.
Let $X$ be the generic fiber of $\Y$, and let $\A$ be the N{\'e}ron model
of the Jacobian of $X$.  Then the connected component 
$\A^0$ of the identity in $\A$ coincides naturally with the
relative Picard scheme $Pic^0_{\Y / R}$.  In particular, every 
automorphism $\psi$ of $\Y$ induces an endomorphism of $\A^0$, and 
if $\psi$ is the identity map on the special fiber of $\Y$ then it 
induces the identity map on the connected component of the closed fiber 
of $\A$.
\end{thm}

\begin{pf}
This is a special case of \cite[Theorem~4, Section 9.5]{BLR}.
\end{pf}

% To see this, let $\modJstar$ be the Neron model of $\Jstar$, 
% and let $\widetilde{u}$ be the induced automorphism of the
% special fiber $\redJstar$.  
% By Lemma \ref{lem:AutX into EndJ} and Proposition \ref{prop:Aut X}, 
% we know that if $u\neq 1$ on 
% $\Xstar$, then $u\neq 1$ on $\Jstar$, and therefore that there is a 
% nontrivial abelian subvariety $J' \subset \Jstar$ on which $u$ acts as
% $-1$.  But then $\widetilde{u} = -1$ on the special fiber 
% $\widetilde{{\mathcal J'}}$,  
% and since the natural map ${\mathcal J'} \to \modJstar$
% is a closed immersion by \cite[Thm.\,7.5.4]{BLR}, it follows that
% $\widetilde{u}\neq 1$ on $\redJstar$.

% Alternatively, using Lemma 2.1, our claim follows from the following 
% general fact: 

\begin{lem}
If $G$ is a semiabelian group scheme over a 
discrete valuation ring $R$, then an endomorphism of $G$ 
inducing the identity map on
the special fiber is necessarily the identity.  
\end{lem}

\begin{pf}

We may assume without loss of generality that $R$ is complete, and we
denote by $\m$ its maximal ideal.

Let $u \in \End G$ be such that $u_s$ (the induced map on the closed
fiber $G_s$ of $G$) is the identity map.  
We want to show that $u=1$ on all of $G$.
It is more convenient to work with the endomorphism $v := u-1$ which 
satisfies $v_s = 0$.  We want to show that $v = 0$.

Fix an auxiliary prime $l$ not equal to the residue characteristic of $R$,
and let $H_n = G[l^n]$ for each integer $n\geq 1$.  
Since $G$ is semiabelian, it follows from 
\cite[Lemma~2, Section~7.3]{BLR} that $H_n$ is a 
quasi-finite, flat, {\'e}tale, and separated group scheme over $R$ for all 
$n$.  Since $R$ is complete, and therefore henselian, we may 
consider, for each $n$, 
the {\em finite part} $(H_n)^{\f}$ of $H_n$,
which is a finite {\'e}tale open and closed subscheme of $H_n$ having the
same closed fiber as $H_n$ (see \cite[Section~7.3, p.~179]{BLR}). 
As $v_s = 0$ on the closed fiber of $(H_n)^{\f}$, and $(H_n)^{\f} / R$ is
finite {\'e}tale, it follows that $v=0$ on $(H_n)^{\f}$.

In particular, $v=0$ on all infinitesimal closed fibers 
of $(H_n)^{\f}$ (i.e., on $(H_n)^{\f} \times_R R/\m^k$ for all integers 
$k\geq 1$).  We claim
that $v=0$ on all infinitesimal closed fibers of $G$.  To see this, note 
that the infinitesimal closed fibers of $(H_n)^{\f}$ and of $H_n$ are the 
same by the construction of $(H_n)^{\f}$.  Also, since $G_s$ is semiabelian,
the collection of subschemes $\{ (H_n)_s \}$ is topologically dense in 
$G_s$.  It follows by \cite[Proposition~11.10.1,Theorem~11.10.9]{EGAIV-3}
that the collection of subschemes induced by $\{ (H_n) \}$ 
is schematically dense in each infinitesimal closed fiber of $G$.  
This proves our claim that $v=0$ on each infinitesimal closed fiber of 
$G/R$. 

It now follows that $v$ induces the zero map on the
formal completion of $G$ along the identity section.  Therefore 
$v=0$ on the generic fiber of $G$, and since $G$ is flat over $R$, 
it follows that $v=0$ on $G$ as desired.
\end{pf}

% THE ARGUMENT WHICH FOLLOWS IS NOT CORRECT AND SHOULD BE DELETED.
% Without loss of generality, we may assume
% that $R$ is complete with residue characteristic $p>0$.  
% Let $T$ be the subgroup of $G$ consisting of all prime-to-$p$ torsion, and
% let $T^\f$ be its finite part.  % Give reference.
% Suppose $u \in \End G$ induces the identity on the 
% special fiber $G_s$.  Then
% $u$ induces the identity on the special fiber $T_s$ of $T$.  
% Since $T^\f$ is finite {\'e}tale over $R$, 
% $u$ induces the identity on $T^\f$.  So it
% suffices to show that $T^\f$ is Zariski-dense in $G$.  For this, let
% $Z$ be the Zariski closure of $T^\f$ in $G$.  If $Z\neq G$, then
% $\dim Z_s < \dim G_s$, which contradicts the fact that $T_s$ is dense in
% the semi-abelian variety $G_s$.

\medskip

Denote by $Y$ the normalization of $\redXstar$, 
which as we have seen is isomorphic to 
$\P^1_{\F_p}$, and let $\{ e_i,e_i' \}$ be
the conjugate pairs of $j$-values of supersingular elliptic curves which
are identified in the map $Y \to \redXstar$
$(i = 1,\ldots,g^*)$.

We conclude from our above discussion that if there is a nontrivial 
element $u\in\Aut \Xstar$, then there is an element $\alpha \in 
\Aut(\P^1_{\F_p})=\PGL(2,\F_p)$ of order 2
% which is defined over $\F_p$, and such that $\alpha$ 
which permutes the set $\{e_i,e_i'\}$ of $2g^*$ supersingular
$j$-invariants which are not $\F_p$-rational.
Now if $g\geq 3$ is an integer, then there is no nontrivial
automorphism of $\P^1$ permuting a general set of $2g$ points.  So if 
$g^* \geq 3$, it is unlikely that an automorphism $\alpha$ as above exists; 
in any case, for each specific value of $p$ one can compute whether or not
such an $\alpha$ exists, since the set $\{e_i,e_i'\}$ can be computed 
explicitly (see Table \ref{tab:j-invariants}).

\medskip

We therefore check by an explicit computation that there is no such element
$\alpha\in\PGL(2,\F_p)$ for the remaining seven cases 
$p=211,223,227,229,389,431,571$.
It follows that $\Aut \Xstar = \{ 1 \}$ for these values of $p$.  This
concludes our proof of Theorem \ref{thm:main}.

\begin{table}[ht]
\caption{Quadratic supersingular $ j $-invariants for the remaining cases}%
\label{tab:j-invariants}
\begin{tabular}{|@{\ }c@{\ }|c|}
\noalign{\hrule height 0.5pt}
$ p $ & Quadratic supersingular $ j $-invariants \\
\noalign{\hrule height 0.5pt}
211 &
$
(j^2{+}162j{+}146)(j^2{+}186j{+}97)(j^2{+}56j{+}23) $ \\[-1.4mm]
& $ \times
(j^2{+}152j{+}88)(j^2{+}184j{+}86)(j^2{+}121j{+}206) $ \\
\noalign{\hrule height 0.1pt}
223 &
$
(j^2{+}87j{+}193)(j^2{+}198j{+}7)(j^2{+}197j{+}1) $ \\[-1.4mm]
& $ \times
(j^2{+}174j{+}49)(j^2{+}137j{+}95)(j^2{+}12j{+}54) $ \\
\noalign{\hrule height 0.1pt}
227 &
$
(j^2{+}168j{+}208)(j^2{+}102j{+}201)(j^2{+}81j{+}63)
(j^2{+}73j{+}163)(j^2{+}223j{+}186) $ \\
\noalign{\hrule height 0.1pt}
229 &
$
(j^2{+}55j{+}175)(j^2{+}14j{+}84)(j^2{+}162j{+}16)
(j^2{+}32j{+}86) $ \\[-1.4mm]
& $ \times
(j^2{+}51j{+}122)(j^2{+}20j{+}2)(j^2{+}63j{+}216) $ \\
\noalign{\hrule height 0.1pt}
389 &
$
(j^2{+}6j{+}326)(j^2{+}233j{+}57)(j^2{+}43j{+}344)
(j^2{+}93j{+}99)(j^2{+}210j{+}377)(j^2{+}240) $ \\[-1.4mm]
& $ \times
(j^2{+}91j{+}374)(j^2{+}27j{+}173)
(j^2{+}293j{+}337)(j^2{+}165j{+}59)(j^2{+}84j{+}36) $ \\
\noalign{\hrule height 0.1pt}
431 &
$
(j^2{+}195j{+}282)(j^2{+}104j{+}148)(j^2{+}301j{+}11)(j^2{+}232j{+}279) $
\\[-1.4mm]
& $ \times
(j^2{+}149j{+}95)(j^2{+}51j{+}138)(j^2{+}10j{+}134)(j^2{+}254j{+}361) $ \\
\noalign{\hrule height 0.1pt}
571 &
$
(j^2{+}513j{+}217)(j^2{+}34j{+}160)(j^2{+}319j{+}194)
(j^2{+}38j{+}51)(j^2{+}284j{+}540) $ \\[-1.4mm]
& $ \times
(j^2{+}12j{+}355)(j^2{+}472j{+}23)(j^2{+}502j{+}522)
(j^2{+}322j{+}545)(j^2{+}25j{+}23) $ \\[-1.4mm]
& $ \times
(j^2{+}307j{+}466)(j^2{+}455j{+}306)(j^2{+}326j{+}7)
(j^2{+}154j{+}51)(j^2{+}30j{+}38) $ \\[-1.4mm]
& $ \times
(j^2{+}222j{+}354)(j^2{+}116j{+}446)(j^2{+}162j{+}68)
(j^2{+}212j{+}138) $ \\
\noalign{\hrule height 0.5pt}
\end{tabular}
\end{table}
%
%
%
%
%
%
%

%
%
\begin{comment}
\vspace{8mm}
{\footnotesize\baselineskip 10pt
Yuji {\sc Hasegawa} \vspace{1mm}\\
{\sc
Department of Mathematical Sciences,
Waseda University \\
3-4-1, Okubo Shinjuku-ku, Tokyo, 169-8555,
Japan}
\\
e-mail: hasegawa@gm.math.waseda.ac.jp
\par
}
\end{comment}
%
%
%
%
\end{document}